\documentclass[11pt,dvips]{article}
\usepackage{amsmath}
\usepackage{amsthm}
\usepackage{graphicx}
\usepackage{bbm,amssymb}
\usepackage{bbold}
\usepackage[hyperfootnotes=false]{hyperref}

\newtheorem{theorem}{Theorem}[section]
\newtheorem{prop}[theorem]{Proposition}
\newtheorem{lemma}[theorem]{Lemma}

\theoremstyle{remark}

\theoremstyle{definition}

\title{The Rotor-Router Model on Regular Trees}
\author{Itamar Landau and Lionel Levine 
\\ University of California, Berkeley}
\date{July 24, 2008}

\DeclareSymbolFont{AMSb}{U}{msb}{m}{n}
\DeclareMathSymbol{\C}{\mathbin}{AMSb}{"43} 
\DeclareMathSymbol{\EE}{\mathbin}{AMSb}{"45} 
\DeclareMathSymbol{\N}{\mathbin}{AMSb}{"4E} 
\DeclareMathSymbol{\PP}{\mathbin}{AMSb}{"50} 
\DeclareMathSymbol{\Q}{\mathbin}{AMSb}{"51} 
\DeclareMathSymbol{\R}{\mathbin}{AMSb}{"52} 
\DeclareMathSymbol{\Z}{\mathbin}{AMSb}{"5A}

\begin{document}
\maketitle
\renewcommand{\thefootnote}{}
\footnote{The second author was supported by a National Science Foundation Graduate Research Fellowship.}
\renewcommand{\thefootnote}{\arabic{footnote}}

\begin{abstract}
The rotor-router model is a deterministic analogue of random walk.  It can be used to define a deterministic growth model analogous to internal DLA.  We show that the set of occupied sites for this model on an infinite regular tree is a perfect ball whenever it can be, provided the initial rotor configuration is acyclic (that is, no two neighboring vertices have rotors pointing to one another).  This is proved by defining the {\it rotor-router group} of a graph, which we show is isomorphic to the sandpile group.  We also address the question of recurrence and transience: We give two rotor configurations on the infinite ternary tree, one for which chips exactly alternate escaping to infinity with returning to the origin, and one for which every chip returns to the origin.  Further, we characterize the possible ``escape sequences" for the ternary tree, that is, binary words $a_1 \ldots a_n$ for which there exists a rotor configuration so that the $k$-th chip escapes to infinity if and only if $a_k=1$.
\end{abstract}

\section{Introduction}

The rotor-router model is a deterministic analogue of random walk, first defined by Priezzhev et al.\ under the name ``Eulerian walkers'' \cite{PDDK} and popularized more recently by Jim Propp \cite{Kleber}.  To define rotor-router walk on a tree $T$, for each vertex of $T$ we choose a cyclic ordering of its neighbors.  Each vertex is assigned a ``rotor'' which points to one of the neighboring vertices.  A chip walks on the vertices of $T$ according to the following rule: when the chip reaches a vertex $v$, the rotor at $v$ rotates to point to the next neighbor in the ordering, and the chip steps in direction of the newly rotated rotor.  In {\it rotor-router aggregation}, we grow a cluster of points in $T$ by repeatedly starting chips at a fixed vertex $o$ and letting them walk until they exit the cluster.  Beginning with $A_1 = \{o\}$, define the cluster $A_n$ inductively by
	\[ A_n = A_{n-1} \cup \{x_n\}, \qquad n >1 \]
where $x_n \in T$ is the endpoint of a rotor-router walk started at $o$ and stopped on first exiting $A_{n-1}$.  We do not change the positions of the rotors when adding a new chip.  Thus the sequence $(A_n)_{n \geq 1}$ depends only on the choice of the initial rotor configuration.  

Recent interest has focused on rotor-router aggregation in the integer lattice $\Z^d$.  Jim Propp noticed from simulations in $\Z^2$ that the shape $A_n$ is extremely close to circular, and asked why this was so \cite{Kleber}.  The spherical shape of $A_n$ in $\Z^d$ is proved in \cite{LP1,LP2}.  Here we prove an analogous result for rotor-router aggregation on the infinite $d$-regular tree.  We say that a rotor configuration is \emph{acyclic} if the rotors form no oriented cycles.  On a tree, this condition is equivalent to forbidding oriented cycles of length~2: there is no pair of neighboring vertices $x,y$ such that both the rotor at $x$ points to $y$ and the rotor at $y$ points to $x$.  As the following result shows, provided we start with an acyclic rotor configuration, the occupied cluster $A_n$ is a perfect ball for suitable values of $n$.

\begin{theorem}
\label{aggregintro}
Let $T$ be the infinite $d$-regular tree, $d \geq 3$, and let
	\[ B_r = \{x \in T \,:\, |x| \leq r\} \]
be the ball of radius $r$ centered at the origin $o \in T$, where $|x|$ is the number of edges in the path from $o$ to $x$.  Write
	\[ b_r = \# B_r = 1 + d \frac{(d-1)^r-1}{d-2}. \]
Let $A_n$ be the region formed by rotor-router aggregation in $T$, starting from $n$ chips at $o$.  If the initial rotor configuration is acyclic, then 
	 \[ A_{b_r} = B_r. \]
\end{theorem}

The proof of Theorem~\ref{aggregintro} uses the {\it sandpile group} of a wired regular tree (that is, a finite regular tree with the leaves collapsed to a single vertex, and an edge added from the root to this vertex), whose structure was found in \cite{sandpiletree}.  In section~2 we define the {\it rotor-router group} of a graph and show that it is isomorphic to the sandpile group.  We then use this isomorphism in section~3 to prove Theorem~\ref{aggregintro}.

Much previous work on the rotor-router model has taken the form of comparing the behavior of rotor-router walk with the expected behavior of random walk.  For example, Cooper and Spencer \cite{CS} show that for any configuration of chips on even lattice sites in $\Z^d$, letting each chip perform rotor-router walk for $n$ steps results in a configuration that differs by only constant error at each point from the expected configuration had the chips performed independent random walks.  In section~4, we continue in this vein by investigating the recurrence and transience of rotor-router walk on trees.  A walk which never returns to the origin visits each vertex only finitely many times, so the positions of the rotors after a walk has escaped to infinity are well-defined.  We construct two ``extremal'' rotor configurations on the infinite ternary tree, one for which walks exactly alternate returning to the origin with escaping to infinity, and one for which every walk returns to the origin.  The latter behavior is something of a surprise: to our knowledge it represents the first example of rotor-router walk behaving fundamentally differently from the expected behavior of random walk.

In between these two extreme cases, a variety of intermediate behaviors are possible.  We say that a binary word $a_1 \ldots a_n$ is an {\it escape sequence} for the infinite ternary tree if there exists an initial rotor configuration on the tree so that the $k$-th chip escapes to infinity if and only if $a_k=1$.  The following result characterizes all possible escape sequences on the ternary tree.

\begin{theorem}
Let $a = a_1 \ldots a_n$ be a binary word. For $j \in \{1,2,3\}$ write $a^{(j)} = a_{j} a_{j+3} a_{j+6} \ldots$. Then $a$ is an escape sequence for some rotor configuration on the infinite ternary tree if and only if for each $j$ and each $k \geq 2$, every subword of $a^{(j)}$ of length $2^k-1$ contains at most $2^{k-1}$ ones.
\end{theorem}

We conclude in section~5 with an open question about the transience of rotor-router walk in $\Z^d$ for $d\geq 3$.

\section{The Rotor-Router Group}

In this section we define the {\it rotor-router group} of a graph and show it is isomorphic to the sandpile group.  The definition of the sandpile group is recalled below.  In the next section we use this isomorphism together with the results of \cite{sandpiletree} to study the rotor-router aggregation model on a regular tree.
%The rotor-router model has also been studied from the viewpoint of self-organized criticality in the physics literature, where it goes by the name ``Eulerian walkers.'' 
The isomorphism between the rotor-router and sandpile groups, Theorem~\ref{groupisom}, is mentioned in the physics literature; see \cite{PDDK, PPS}.  To our knowledge the details of the proof are not written down anywhere.  While our main focus is on the tree, the isomorphism is just as easily proved for general graphs, and it seems to us worthwhile to record the general proof here.

Let $G$ be a strongly connected finite directed graph, which may have multiple edges but not loops.  Fix a vertex $s$ in $G$ and call it the sink.  To define rotor-router walk on $G$, for each vertex $x\neq s$ we fix a cyclic ordering of the edges emanating from $x$.  A \emph{rotor configuration} $T$ on $G$ assigns to each non-sink vertex $x$ an edge $T(x)$ emanating from $x$.    Each step of the walk then consists of two parts: If the chip is located at $x$, we first increment the rotor $T(x)$ to the next edge in the ordering of the edges emanating from~$x$, and then move the chip along this new edge.
Given a rotor configuration~$T$, write $e_x(T)$ for the rotor configuration resulting from starting a chip at $x$ and letting it walk according to the rotor-router rule until it reaches the sink.  (Note that if the chip visits a vertex infinitely often, it visits all of its outbound neighbors infinitely often; since $G$ is strongly connected, the chip eventually reaches the sink.)  
  
The set of edges $\{T(x)\}_{x\neq s}$ in a rotor configuration forms a spanning subgraph of $G$ in which every vertex except the sink has out-degree one.  If this subgraph contains no directed cycles (equivalently, no cycles), we call it an \emph{oriented spanning tree} of $G$.  Write $Rec(G)$ for the set of oriented spanning trees of $G$.  Note that as we have defined them, oriented spanning trees are always rooted at the sink (i.e., all paths in the tree lead to the sink).

\begin{lemma}
If $T \in Rec(G)$, then $e_x(T) \in Rec(G)$.
\end{lemma}

\begin{proof}
Let $Y$ be any collection of vertices of $G$.  If the chip started at $x$ reaches the sink without ever visiting $Y$, then the rotors at vertices in $Y$ point the same way in $e_x(T)$ as they do in $T$, so they do not form an oriented cycle.  If the chip does visit $Y$, let $y\in Y$ be the last vertex it visits.  Then either $y=s$, or the rotor at $y$ points to a vertex not in $Y$; in either case, the rotors at vertices in $Y$ do not form an oriented cycle.
\end{proof}

We will need slightly more refined information about the intermediate states that occur before the chip falls into the sink. These states may contain oriented cycles, but only of a very restricted form.  For a vertex $x$ we write $Cyc_x(G)$ for the set of rotor configurations $U$ such that
	\begin{itemize}
	\item[(i)] $U$ contains an oriented cycle; and
	\item[(ii)] If the rotor $U(x)$ is deleted, the resulting subgraph contains no oriented cycles.
	\end{itemize}

\begin{lemma}
\label{cyclecriterion}
Starting from a rotor configuration $T_0 \in Rec(G)$ with a chip at $x_0$, let $T_k$ and $x_k$ be the rotor configuration and chip location after $k$ steps of rotor-router walk.  Then
\begin{itemize}
\item[(i)] If $T_k \notin Rec(G)$, then $T_k \in Cyc_{x_k}(G)$.
\item[(ii)] If $T_k \in Rec(G)$, then $x_k \notin \{x_0,\ldots,x_{k-1}\}$.
%\item[(i)] If $x_k \in \{x_0, \ldots, x_{k-1}\}$, then $T_k$ has an oriented cycle containing $x_k$, and if the rotor at $x_k$ is deleted, the resulting graph has no oriented cycles.
%\item[(ii)] If $x_k \notin \{x_0, \ldots, x_{k-1}\}$, then $T_k \in Rec(G)$.
\end{itemize}
\end{lemma}

\begin{proof}
(i) It suffices to show that any oriented cycle in~$T_k$ contains~$x_k$.  Let~$Y$ be any set of vertices of~$G$ not containing~$x_k$.  If~$Y$ is disjoint from $\{x_0, \ldots, x_{k-1}\}$, then the rotors at vertices in~$Y$ point the same way in~$T_k$ as they do in~$T_0$, so they do not form an oriented cycle.  Otherwise, let $y\in Y$ be the vertex visited latest before time~$k$.  The rotor~$T_k(y)$ points to a vertex not in~$Y$, so the rotors at vertices in~$Y$ do not form an oriented cycle.

(ii) Suppose $x_k \in \{x_0,\ldots,x_{k-1}\}$.  Let~$y_0=x_k$, and 
%define $y_1,y_2,\ldots$ by $T_k(y_i)=(y_i,y_{i+1})$.
for $i=0,1,\ldots$ let~$y_{i+1}$ be the target of the rotor~$T_k(y_i)$.  Then the last exit from~$x_k$ before time~$k$ was to~$y_1$, and by induction if $y_1, \ldots, y_{i-1}$ are different from~$x_k$, then~$y_{i-1}$ was visited before time~$k$, and the last exit from~$y_{i-1}$ before time~$k$ was to~$y_i$.  It follows that~$y_i=x_k$ for some~$i\geq 1$, and hence~$T_k$ contains an oriented cycle.
\end{proof}

\begin{lemma}
If $T_1, T_2 \in Rec(G)$ and $e_x(T_1)=e_x(T_2)$, then $T_1=T_2$.
\end{lemma}

\begin{proof}
We will show that any $T\in Rec(G)$ can be recovered from $e_x(T)$ by reversing one rotor step at a time.
Given rotor configurations $U,U'$ and vertices $y,y'$, we say that $(U',y')$ is a predecessor of $(U,y)$ if a chip at $y'$ with rotors configured according to $U'$ would move to $y$ in a single step with resulting rotors configured according to $U$.  Given $U$ and $y$, for each neighbor $z$ of $y$ whose rotor $U(z)$ points to $y$, there is a unique predecessor of the form $(U',z)$, which we will denote $P_z(U,y)$.
% fine for multigraphs too

Suppose $(U,y)$ is an intermediate state in the evolution from $T$ to $e_x(T)$.  If $U \notin Rec(G)$, then by case (i) of Lemma~\ref{cyclecriterion} there is a cycle of rotors
$y \rightarrow y_1 \rightarrow y_2 \rightarrow \ldots \rightarrow y_n \rightarrow y$ in $U$.
%$U(y)=(y,y_1), U(y_1)=(y_1,y_2), \ldots, U(y_n)=(y_n,y)$.  
If $z$ is a vertex different from $y_n$ whose rotor $U(z)$ points to $y$, then $z$ is not in this cycle, so the predecessor $P_z(U,y)$ has a cycle disjoint from its chip location.  Thus $P_z(U,y)$ does not belong to $Rec(G)$ or to $Cyc_z(G)$, so by Lemma~\ref{cyclecriterion} it cannot be an intermediate state in the evolution from $T$ to $e_x(T)$.  The state immediately preceding $(U,y)$ in the evolution from $T$ to $e_x(T)$ must therefore be $P_{y_n}(U,y)$.

Now suppose $U \in Rec(G)$.  By case (ii) of Lemma~\ref{cyclecriterion}, $U$ is the rotor configuration when $y$ is first visited.  If $y=x$, then $U=T$.  Otherwise, let $x=x_0 \rightarrow x_1 \rightarrow \ldots \rightarrow x_k=s$ be the path in $U$ from $x$ to the sink.  Then the last exit from $x$ before visiting $y$ was to $x_1$.  By induction, if $x_1, \ldots, x_{j-1}$ are different from $y$, then $x_{j-1}$ was visited before $y$ and the last exit from $x_{j-1}$ before visiting $y$ was to $x_j$.  It follows that $x_j=y$ for some $j\geq 1$, and the state immediately preceding $(U,y)$ must be $P_{x_{j-1}}(U,y)$.
\end{proof}

Thus for any vertex $x$ of $G$, the operation $e_x$ of adding a chip at $x$ and routing it to the sink acts invertibly on the set of states $Rec(G)$ whose rotors form oriented spanning trees rooted at the sink.  It is for this reason that we call these states recurrent.  We define the {\it rotor-router group} $RR(G)$ as the subgroup of the permutation group of $Rec(G)$ generated by $\{e_x\}_{x\neq s}$.  For any two vertices $x$ and $y$, the operators $e_x$ and $e_y$ commute; this commutativity is proved in \cite{DF} for a broad class of models encompassing both the abelian sandpile and the rotor-router.
%Note that if there are two (indistinguishable) chips on $G$ and each takes a single step according to the rotor-router rule, the resulting rotor configuration does not depend on the order of the two steps.  A somewhat more subtle argument shows that the operators $e_x$ commute
Hence the group $RR(G)$ is abelian.

\begin{lemma}
\label{transitivity}
$RR(G)$ acts transitively on $Rec(G)$.
%For any $T_1, T_2 \in Rec(G)$ there exists $u$ such that $T_2 = \phi(u) T_1$.
\end{lemma}

\begin{proof}
Given $T_1, T_2 \in Rec(G)$, for each vertex $x\neq s$ let $u(x)$ be the number of rotor turns needed to get from $T_1(x)$ to $T_2(x)$.  Let $v(x)$ be the number of chips ending up at~$x$ if $u(y)$ chips start at each vertex~$y$, with rotors starting in configuration~$T_1$, and each chip takes a single step.  After each chip has taken a single step, the rotors are in configuration $T_2$, hence
		\[ \left( \prod_{x \neq s} e_x^{u(x)} \right) T_1 =  \left( \prod_{x \neq s} e_x^{v(x)} 
\right) T_2. \]

%		\[ \left( \sum_{x \in V(G)} u(x)e_x \right) T_1 =  \left( \sum_{x \in V(G)} v(x)e_x \right) T_2. \]
Letting $g = \prod_{x \neq s} e_x^{u(x)-v(x)}$ we obtain $T_2 = gT_1$.
\end{proof}

Given vertices $x$ and $y$, write $d_{xy}$ for the number of edges in $G$ from $x$ to $y$, and write
	\[ d_x = \sum_y d_{xy} \]
for the outdegree of $x$.

\begin{theorem}
\label{groupisom}
Let $G$ be a strongly connected finite directed graph without loops, let $RR(G)$ be its rotor-router group, and $SP(G)$ its sandpile group.  Then $RR(G) \simeq SP(G)$.
\end{theorem}

\begin{proof}
Let $V$ be the vertex set of $G$.  The sandpile group of $G$ \cite{CR,Dhar} 
% CR only consider undirected graphs!
is the quotient
	\[ SP(G) = \Z^V \big/ (s,\Delta_x)_{x\in V} \]
where $s \in V$ is the sink and
	\[ \Delta_x = \sum_{y\in V} d_{xy} y -  d_x x. \]
Define $\phi : \Z^V \rightarrow RR(G)$ by
	\[ \phi \left(\sum_{x\in V} u_x x\right) = \prod_{x \in V} e_x^{u_x}. \]
Starting with $d_x$ chips at a vertex $x$ and letting each chip take one rotor-router step results in $d_{xy}$ chips at each vertex $y$, with the rotors unchanged, hence
	\[ e_x^{d_x} = \prod_{y\in V} e_y^{d_{xy}}. \]
Thus $\phi(\Delta_x)=Id$.  Since also $\phi(s)=e_s=Id$, the map $\phi$ descends to a map $\bar{\phi}:SP(G) \rightarrow RR(G)$.  This latter map is surjective since $\phi$ is surjective; to show that $\bar{\phi}$ is injective, by Lemma~\ref{transitivity} we have
	\[ \# RR(G) \geq \# Rec(G) = \# SP(G), \]
where the equality on the right is the matrix-tree theorem \cite[5.6.8]{Stanley}. 
\end{proof}

\section{Aggregation on the Tree}

Fix $d\geq 3$, and let $T$ be the infinite $d$-regular tree.  Fix an origin vertex~$o$ in~$T$.  In {\it rotor-router aggregation}, we grow a cluster of points in~$T$ by repeatedly starting chips at the origin and letting them walk until they exit the cluster.  Beginning with $A_1 = \{o\}$, define the cluster~$A_n$ inductively by
	\[ A_n = A_{n-1} \cup \{x_n\}, \qquad n >1. \]
where $x_n \in T$ is the endpoint of a rotor-router walk started at~$o$ and stopped on first exiting~$A_{n-1}$.  We do not change the positions of the rotors when adding a new chip.  In this section we use the group isomorphism in Theorem~\ref{groupisom} to show that $A_n$ is a perfect ball for suitable values of $n$ (Theorem~\ref{treecirc}).  

A function $H$ on the vertices of a directed graph $G$ is {\it harmonic} at a vertex $x$ if
	\[ d_x H(x) = \sum_{y \in V} d_{xy} H(y), \]
where $d_{xy}$ is the number of edges from~$x$ to~$y$, and~$d_x$ is the outdegree of~$x$.
%	\[ H(x) = \frac{1}{\text{outdeg}(x)} \sum_{y \leftarrow x} H(y) \]

\begin{lemma}
\label{HPinvariant}
Let $G=(V,E)$ be a finite directed graph without loops.
Suppose chips on $G$ can be moved by a sequence of rotor-router steps, starting with $u(x)$ chips at each vertex $x$ and ending with $v(x)$ chips at each vertex $x$, in such a way that the initial and final rotor configurations are the same.  If $H$ is a function on $V$ that is harmonic at all vertices which emitted chips, then
	\[ \sum_{x \in V} H(x) u(x) = \sum_{x \in V} H(x) v(x). \]
\end{lemma}

\begin{proof}
Let $u=u_0, u_1, \ldots, u_k=v$ be the intermediate configurations.  If $u_{i+1}$ is obtained from $u_i$ by routing a chip from $x_i$ to $y_i$, then
	\begin{equation}\label{stepbystep} \sum_{x \in V} H(x) (u(x)-v(x)) = \sum_i H(x_i)-H(y_i). \end{equation}
If the initial and final rotor configurations are the same, then each rotor makes an integer number of full turns, so the sum in (\ref{stepbystep}) can be written 
	\[ \sum_i H(x_i)-H(y_i) = \sum_{x \in V} N(x) \sum_{y \in V} d_{xy} (H(x)-H(y)) \]
where $N(x)\in \Z_{\geq 0}$ is the number of full turns performed by the rotor at~$x$.  By the harmonicity of $H$, the inner sum on the right vanishes whenever $N(x)>0$.
\end{proof}

Next we describe our choice of graph $G$ and harmonic function $H$.
By the \emph{$d$-regular tree of height~$n$} we will mean the finite rooted tree in which each non-leaf vertex has $d-1$ children, and the path from each leaf to the root has $n-1$ edges.  We denote this tree by~$T_n$.  Let $\hat{T}_n$ be the graph obtained from $T_n$ by adding a single additional leaf $o$ whose parent is the root $r$ of $T_n$.  This is an undirected graph; when applying the results above, which are phrased in terms of directed graphs for maximum generality, we think of it as \emph{bidirected}: each edge is replaced by a pair of directed edges pointing in opposite directions.

Denote by $(X_t)_{t\geq 0}$ the simple random walk on $\hat{T}_n$, and let $\tau \geq 0$ be the first hitting time of the set of leaves.  Fix a leaf $z\neq o$, and let
	\begin{equation} \label{ourharmonicfunction} H(x) = \PP_x(X_\tau = z) \end{equation}
be the probability that random walk started at $x$ and stopped at time $\tau$ stops at $z$.  This function is harmonic at all non-leaf vertices.

We briefly recall the well-known martingale argument from gambler's ruin used to find the value of $H(r)$.  The process
	\[ M_t = a^{-|X_t|} \]
is a martingale, where $a=d-1$ and $|x|$ denotes the number of edges in the path from $o$ to $x$.  Since $M_t$ has bounded increments and $\EE_r \tau < \infty$, we obtain from optional stopping
	\[ a^{-1} = \EE_r M_0 = \EE_r M_\tau = p + (1-p) a^{-n} \]
where $p = \PP_r(X_\tau = o)$.  Solving for $p$ we obtain
	\begin{equation} \label{gamblersruin} \PP_r(X_\tau=o) = \frac{a^{n-1}-1}{a^n-1}. \end{equation}
In the event that the walk stops at a leaf $z\neq o$, by symmetry it is equally likely to stop at any such leaf.  Since there are $a^{n-1}$ such leaves, we obtain from (\ref{gamblersruin})
	\begin{equation} \label{hittingprobofleaf} H(r) = \frac{1 - \PP_r(X_\tau = o)}{a^{n-1}}
										= \frac{a-1}{a^n-1}. \end{equation}

The \emph{wired $d$-regular tree} of height $n$ is the graph~$\bar{T}_n$ obtained from~$\hat{T}_n$ by collapsing all the leaves to a single vertex~$s$, the {\it sink}.  We do not collapse edges; thus each neighbor of the sink except for $r$ has $a=d-1$ edges to the sink.  The proof of Theorem~\ref{treecirc} will use the following fact about the sandpile group of the wired regular tree.

\begin{lemma}
\label{degreeofroot}
The root $r$ of $\bar{T}_n$ has order $\frac{a^n-1}{a-1}$ in the sandpile group $SP\big(\bar{T}_n\big)$.  
\end{lemma}

\begin{proof}
See \cite{sandpiletree}, Proposition~4.2.
\end{proof}

The next lemma concerns rotor-router walk on $\hat{T}_n$ stopped on hitting the leaves.  The leaves play the role of sinks, and the dynamics are the same as for rotor-router walk on the wired tree $\bar{T}_n$.  However, we are interested in counting how many chips stop at each leaf, which is why we preserve the distinction between $\hat{T}_n$ and $\bar{T}_n$.  Since the rotors at the leaves play no role, we view our rotor configuration as living on $\bar{T}_n$.  Such a configuration is \emph{acyclic} if no two neighboring vertices have rotors pointing to one another; in the notation of the previous section, the acyclic configurations are precisely those in $Rec(\bar{T}_n)$.

\begin{lemma}
\label{exitmeasure}
Let $a=d-1$.  Given an acyclic rotor configuration on $\bar{T}_n$, starting with $\frac{a^n-1}{a-1}$ chips at the root $r$ of $\hat{T}_n$, and stopping each chip when it reaches a leaf, exactly one chip stops at each leaf $z\neq o$, and the remaining $\frac{a^{n-1}-1}{a-1}$ chips stop at $o$.  Moreover, the starting and ending rotor configurations are identical.
\end{lemma}

\begin{proof}
By Theorem~\ref{groupisom} and Lemma~\ref{degreeofroot}, the element $e_r \in RR(\bar{T}_n)$ has order $m=\frac{a^n-1}{a-1}$, so $e_r^m$ is the identity permutation of $Rec(\bar{T}_n)$, hence the starting and ending rotor configurations are identical.  Fix a leaf $z\neq o$ of $\hat{T}_n$ and let~$H$ be the function on vertices of $\hat{T}_n$ given by (\ref{ourharmonicfunction}).  Since $H$ is harmonic on the non-leaf vertices, by Lemma~\ref{HPinvariant} and (\ref{hittingprobofleaf}), the number of chips stopping at $z$ is 
	\[  \sum H(x) v(x) = \sum H(x) u(x) = \frac{a^n-1}{a-1} H(r) = 1. \]
Since there are $a^{n-1}$ leaves $z\neq o$, the remaining $\frac{a^n-1}{a-1}-a^{n-1} = \frac{a^{n-1}-1}{a-1}$ chips stop at o.
\end{proof}

The \emph{principal branches} of the infinite $d$-regular tree $T$ are the $d$ subtrees rooted at the neighbors of the origin.  The ball of radius $\rho$ centered at the origin in $o \in T$ is
	\[ B_\rho = \{x \in T \,:\, |x| \leq \rho \} \]
where $|x|$ is the number of edges in the path from $o$ to $x$.  Write
	\[ b_\rho = \# B_\rho = 1 + (a+1) \frac{a^\rho-1}{a-1}. \]
As the following result shows, provided we start with an acyclic configuration of rotors, the rotor-router aggregation cluster $A_n$ is a perfect ball at those times when an appropriate number of chips have aggregated.  It follows that at all other times, the cluster is as close as possible to a ball: if $b_\rho<n<b_{\rho+1}$ then $B_\rho \subset A_n \subset B_{\rho+1}$.

\begin{theorem}
\label{treecirc}
Let $A_n$ be the region formed by rotor-router aggregation on the infinite $d$-regular tree, starting from $n$ chips at the origin.  If the initial rotor configuration is acyclic, then $A_{b_\rho} = B_\rho$ for all $\rho \geq 0$.
%	\begin{itemize}
%	\item[(i)] If $n=b_r$, then $A_n=B_r$.
%, and the rotors at time $n$ are in their initial state. 
%	\item[(ii)] If $b_r<n<b_{r+1}$ then $B_r \subset A_n \subset B_{r+1}$.
%	\end{itemize}
\end{theorem}

\begin{proof}
Define a modified aggregation process $A'_n$ as follows.  Stop the $n$-th chip when it either exits the occupied cluster $A'_{n-1}$ or returns to $o$, and let
	\[ A'_n = A'_{n-1} \cup \{x'_n\} \]
where $x'_n$ is the point where the $n$-th chip stops.  By relabeling the chips, this yields a time change of the original process, i.e.\ $A'_n = A_{f(n)}$ for some sequence $f(1), f(2), \ldots$.  Thus it suffices to show $A'_{c_\rho}=B_\rho$ for some sequence $c_1, c_2, \ldots$.  We will show by induction on $\rho$ that this is the case for
	\[ c_\rho = 1 + (a+1) \sum_{t=1}^{\rho} \frac{a^t-1}{a-1}, \]
and that after $c_\rho$ chips have stopped, the rotors are in their initial state.  For the base case $\rho=1$, we have $c_1 = a+2=d+1$.  The first chip stops at $o$, and the next $d$ stop at each of the neighbors of $o$, so $A'_{d+1}=B_1$.  Since the rotor at $o$ has performed one full turn, it is back in its initial state.

Assume now that $A'_{c_{\rho-1}} = B_{\rho-1}$ and that the rotors are in their initial acyclic state.  Starting with $c_\rho - c_{\rho-1}$ chips at $o$, let each chip in turn perform rotor-router walk until either returning to $o$ or exiting the ball $B_{\rho-1}$.  Then each chip is confined to a single principal branch of the tree, and each branch receives~$\frac{a^\rho-1}{a-1}$ chips.
By Lemma~\ref{exitmeasure}, exactly one chip will stop at each leaf $z \in B_{\rho}-B_{\rho-1}$, and the remainder will stop at $o$.  Thus $A'_{c_\rho} = B_\rho$.  Moreover, by Lemma~\ref{exitmeasure}, once all chips have stopped, the rotors are once again in their initial state, completing the inductive step.
\end{proof}

\section{Recurrence and Transience}

\begin{figure}
	\centering
		\includegraphics[scale=.7]{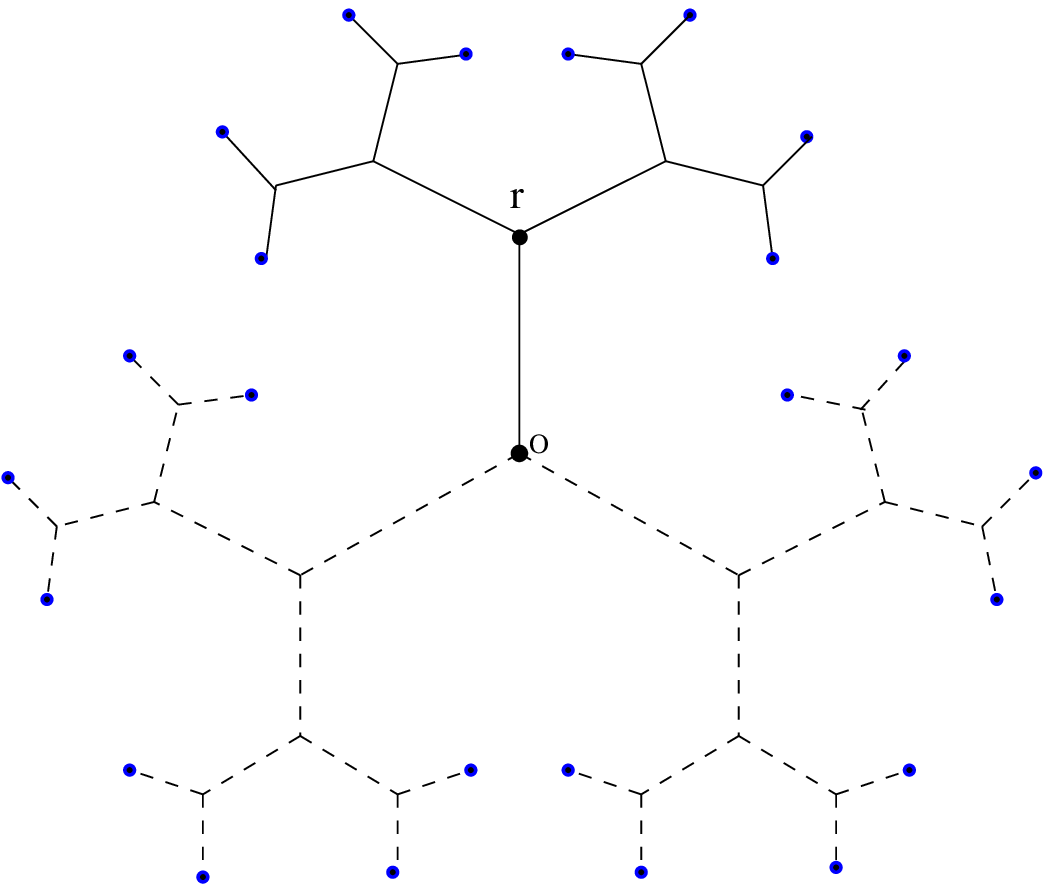}
		\\ ~ \\ ~ \\
		\includegraphics[scale=.7]{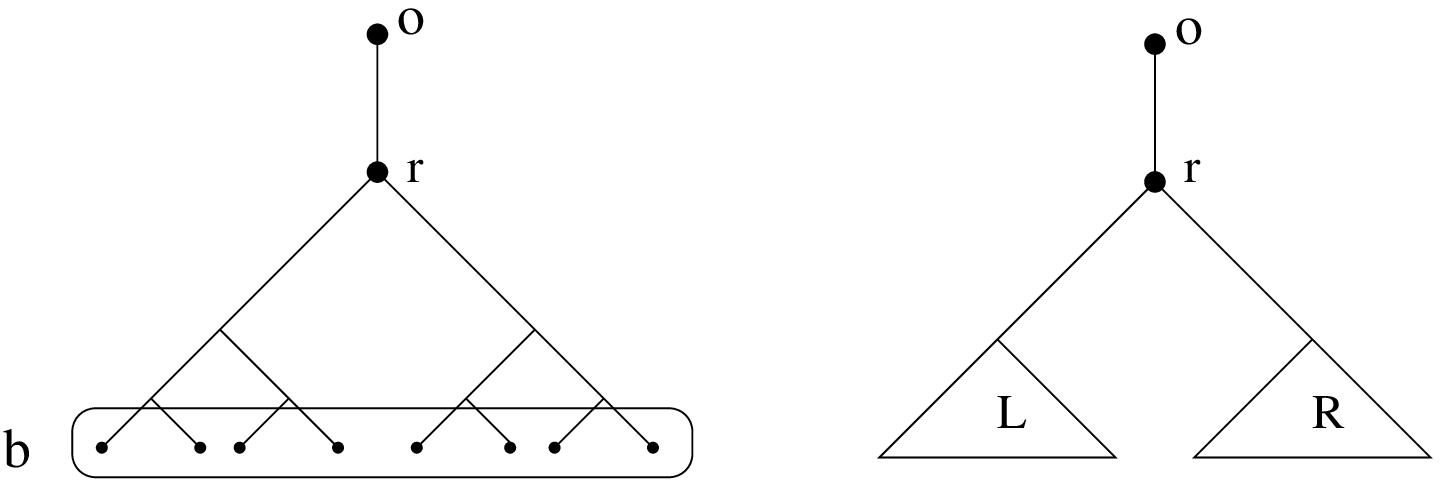}
	\caption{The ball $B_n$ in the regular ternary tree (top), the branch $Y_n$ (left), and its sub-branches $L$ and $R$.}
	\label{fig:Tree1}
\end{figure}

In this section we explore questions of recurrence and transience for the rotor-router walk on regular trees. We aim to study to what extent the rotor-router walk behaves as a deterministic analogue of random walk.  We find that the behavior depends quite dramatically on the initial configuration of rotors.

A chip performing rotor-router walk starting at the origin $o$ in the infinite $d$-regular tree either returns to the origin or escapes to infinity within a single principal branch of the tree, leaving the rotors in the other branches unchanged.
 Therefore, as shown in Figure~1, we focus on a single branch~$Y_n$ of the ball~$B_n$ in the $d$-regular tree.  In the notation of the previous section,~$Y_n$ is the graph obtained from $\hat{T}_n$ by collapsing all the leaves except for $o$ to a single vertex, which we label $b$ for boundary.  Starting chips at the root~$r$ of~$Y_n$, and stopping them either when they reach~$b$ or return to~$o$, we will compare the hitting rates of~$o$ and~$b$ for rotor-router walk with the expected hitting rates for random walk.  
 
To each rotor direction we associate an index from $\{1,\ldots,d\}$, with direction $d$ corresponding to a rotor pointing to the parent vertex.  Rotors cycle through the indices in order.  In the ternary tree ($d=3$) we will sometimes refer to the three rotor directions as left (direction $1$), right (direction $2$) and up (direction $3$).

\begin{lemma}
\label{ternarycase}
Suppose $d=3$.
If all rotors in $Y_n$ initially point in direction~$1$, then the first $2^{n}-1$ chips started at~$r$ alternate, the first stopping at~$b$, the next stopping at~$o$, the next at~$b$, and so on.  After this sequence of $2^n-1$ walks, all rotors again point in direction~$1$.
\end{lemma}

\begin{proof}
Induct on $n$.  In the base case $n=2$, there is only one rotor, which sends the first chip in direction $2$ to $b$, the next chip up in direction $3$ to $o$, and the third chip in direction $1$ to $b$, at which point the rotor is again in its initial state.

Now suppose that the lemma holds for $Y_{n-1}$.  Let $L$ and $R$ be the two principal branches of $Y_n$.  We think of $L$ and $R$ as each having a rotor that points either to $b$ or back up to $r$. The initial state of these rotors is pointing to $r$. The first
chip is sent from the root to $R$, which by induction sends it to $b$. Note that the root rotor is now pointing towards $R$, the $R$-rotor is pointing to $b$, and the $L$-rotor is pointing to $r$ (Figure~\ref{fig:FourSteps}a). We now observe that the next four chips form a
pattern that will be repeated. The second chip is sent directly to $o$ (Figure~\ref{fig:FourSteps}b) and
the third chip is sent to $L$ which sends it to $b$ (Figure~\ref{fig:FourSteps}c). The fourth chip
is sent to $R$, but by induction this chip is returned and then it is sent
to $o$ (Figure~\ref{fig:FourSteps}d). Finally, the fifth chip is sent to $L$, returned, sent to $R$, and
through to~$b$ (Figure~\ref{fig:FourSteps}e). Note that the root rotor is now again pointing towards $R$,
the $R$-rotor is again pointing to~$b$, and the $L$-rotor is again pointing to $r$. In this cycle of four chips, the two branches $R$ and $L$ see two chips apiece. This cycle repeats $2^{ n - 2 } - 1$
times, and each subtree sees $2^{ n - 1 } - 2$ chips.

\begin{figure}
	\centering
		\includegraphics[width=\textwidth]{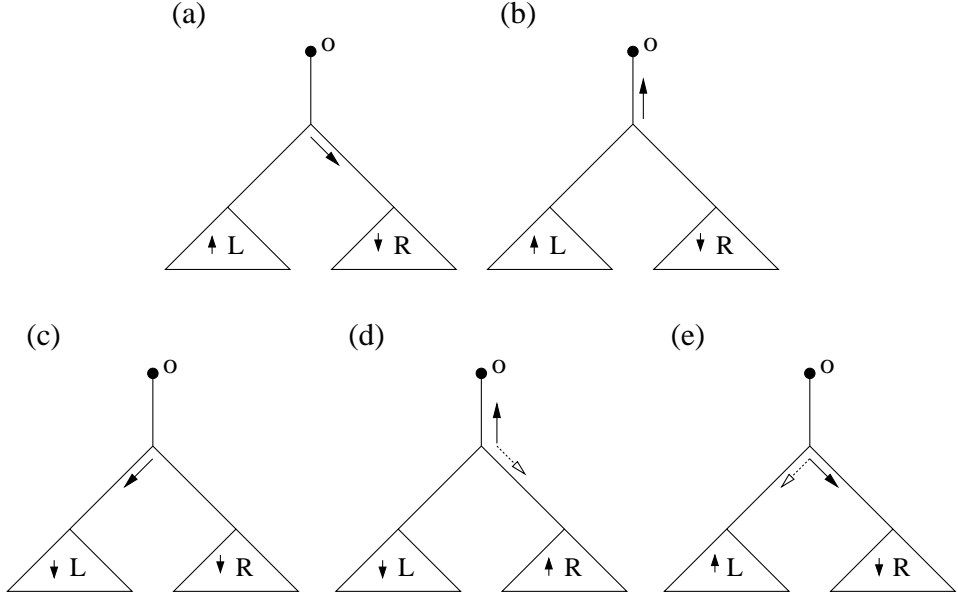}
	\caption{The four-chip cycle, which begins after the first chip has been routed to $b$.}
	\label{fig:FourSteps}
\end{figure}

Recall that the first chip was sent to $R$, so $R$ it has seen a total of $2^{ n - 1 }
- 1$ chips. By induction, all the rotors in $R$ are in their initial
configuration. We have sent a total of $2^n - 3$ chips. The next chip
is sent to $o$, and the last to $L$, which sends it to $b$. Now $L$ has seen
$2^{ n - 1 } - 1$ chips, so by induction all of its rotors are in their
initial configuration. The root rotor is pointing towards $L$, its initial
configuration. We have sent a total of $2^n - 1$ chips, alternating
between $b$ and $o$, and all of the rotors of $Y_n$ are in the initial
configuration, so the inductive step is complete.
\end{proof}

We remark that the obvious generalization of Lemma~\ref{ternarycase} to trees of degree $d > 3$ fails; indeed, we do not know of a starting rotor configuration on trees of higher degree which results in a single chip stopping at $o$ alternating with a string of $d-1$ chips stopping at $b$.

Consider now the case of the infinite ternary tree $T$.  A chip performing rotor-router walk started at the origin $o\in T$ must either return to $o$ or escape to infinity visiting each vertex only finitely many times. 
%For suppose a chip executes an infinite cycle over a finite subset of the tree. Then there must be a vertex which the chip hits infinitely many times. But then the chip will also hit each of the neighboring verteces infinitely many times, and by induction it must hit every vertex infinitely many often. In particular, it must hit and stop at $o$. 
Thus the state of the rotors after a chip has escaped to infinity is well-defined.  We can therefore run a sequence of~$m$ rotor-router walks and count the number~$R(m)$ that return to the origin.  The following result shows that there is an initial rotor configuration on the tree for which the rotor-router walk behaves as an exact quasirandom analogue to the random walk, in which chips exactly alternate returning to the origin with escaping to infinity.

\begin{prop}
\label{ternarythm}
Let $T$ be the infinite ternary tree, with principal branches labeled $Y^{(1)}$, $Y^{(2)}$, and $Y^{(3)}$ in correspondence with the direction indexing of the rotor at the origin. Set the rotors along the rightmost path to infinity in $Y^{(3)}$ initially pointing in direction~$2$, and all remaining rotors initially pointing in direction~$1$. Let $E(m)$ be the expected number of chips that return to the origin if $m$ chips perform independent random walks on~$T$. Let $R(m)$ be the number of chips that return to the origin if $m$ chips sequentially perform rotor-router walks on~$T$. Then $|E(m) - R(m)| \leq \frac{1}{2}$ for all $m$.
\end{prop}

\begin{proof}
Lemma~\ref{ternarycase} implies that for the branches $Y^{(1)}$ and $Y^{(2)}$, the chips sent to a given branch alternate indefinitely with the first escaping to infinity, the next returning to $o$, and so on.  Likewise, chips sent to $Y^{(3)}$ will alternate indefinitely with the first returning to $o$, the next escaping to infinity, and so on. Since chips on the full tree $T$ are routed cyclically through the branches beginning with $Y^{(2)}$, we see that the chips too will alternate indefinitely between escaping to infinity and returning to the origin, with the first escaping to infinity. Thus $R(m) = \left\lfloor \frac{m}{2} \right\rfloor$. Taking $n \rightarrow \infty$ in (\ref{gamblersruin}) we obtain $E(m) = \frac{m}{2}$, and the result follows.
\end{proof}

\begin{lemma}
\label{finitetreesreturn}
For any $d\geq 3$, if all rotors in~$Y_n$ initially point in direction $d-1$, then the first $n-1$ chips started at~$r$ all hit~$o$ before hitting~$b$.  After these $n-1$ chips have stopped at~$o$, the final rotors all point in direction~$d$.
\end{lemma}

\begin{proof}
Induct on $n$. In the base case $n=2$, the first chip steps directly from $r$ to $o$, leaving the single rotor pointing in direction $d$.
Now suppose the lemma holds for $Y_{n-1}$. Let $Z_1, \ldots, Z_{d-1}$ be the principal branches of $Y_n$.  The first chip placed at~$r$ is sent directly to~$o$. By the inductive hypothesis, the first $n-2$ chips that are sent to each branch $Z_i$ are returned to~$r$ before hitting~$b$.  Thus each of the next $n-2$ chips started at~$r$ is sent to $Z_1$, returned to~$r$, sent to $Z_2$, and so on until it is sent to $Z_{d-1}$, returned to~$r$ and then routed to~$o$.  The root rotor now points in direction~$d$, and since each branch~$Z_i$ received exactly $n-2$ chips, its final rotors all point in direction~$d$ by the inductive hypothesis.
\end{proof}

Our next result shows that, perhaps surprisingly, the initial rotors can be set up so as to make rotor-router walk on the $d$-regular tree recurrent.

\begin{prop}
\label{infinitetreesreturn}
On the infinite $d$-regular tree $T$, if all rotors initially point in direction $d-1$, then every chip in an infinite succession of chips started at the origin eventually returns to the origin.
\end{prop}

\begin{proof}
By Lemma~\ref{finitetreesreturn}, for each $n$, the $n$-th chip sent to each principal branch~$Y$ returns to the origin before hitting height $n+1$ of~$T$.
\end{proof}

Note also that if all the rotors in the first $n-1$ levels of~$T$ initially point in direction $d-1$, and all remaining rotors initially point in direction~$d$, then after $n-1$ chips have been sent to a given branch $Y$ and returned to the origin, by Lemma~\ref{finitetreesreturn} all rotors in $Y$ point in direction~$d$, so the next chip sent to~$Y$ escapes to infinity.
%needed for the proof of Lemma~\ref{branch-escape}

We continue our exploration of recurrence and transience on the infinite ternary tree~$T$, allowing now for arbitrary rotor configurations. We focus on a single principal branch~$Y$ of the infinite tree, rooted at a neighbor~$r$ of the origin $o\in T$.  We include the edge~$(o,r)$ in~$Y$, so that~$r$ has degree~$d$ in~$Y$, and~$o$ has degree one.  Thus each chip started at the origin will move to $r$ on its first step.  Given a rotor configuration on~$Y$, we define the {\it escape sequence} for the first~$n$ chips to be the binary word $a = a_{1}\ldots a_{n}$, where for each~$j$, 
\[a_{j} = \begin{cases} 0, & \text{if the $j^{th}$ chip returns to the origin;} \\
				    1, & \text{if the $j^{th}$ chip escapes to infinity.} \end{cases} \]
As noted previously, a chip cannot stay within a finite height indefinitely without returning to the origin, so $a$ is well-defined.

We define a map $\psi$ associating to an escape sequence $a = a_{1} \ldots a_{n}$ a pair of shorter sequences. First, we rewrite $a$ as the concatenation of subwords $b_{1} \cdots b_{m}$ where each $b_{j} \in \{0,10,110\}$. Since at least one of any three consecutive chips entering~$Y$ is routed back to the origin by the rotor at the root~$r$ of~$Y$, at most two of any three consecutive letters in an escape sequence~$a$ can be ones.  Therefore, any escape sequence can be factored in this way up to the possible concatenation of an extra $0$. Now we define $\psi(a) = (c,d)$ by
\begin{equation} \label{psidef} (c_{j},d_{j}) = \begin{cases} (0,0), & \text{if} \; b_{j} = 0\\
																(1,1), & \text{if} \; b_{j} = 110\\
																(0,1), & \text{if} \; b_{j} = 10 \; \text{and} \; \#\{i<j|b_{i}=10\} \; \text{is odd}\\
																(1,0), & \text{if} \; b_{j} = 10 \; \text{and} \; \#\{i<j|b_{i}=10\} \; \text{is even.} \end{cases} 
\end{equation}

In the other direction, given a pair of binary words $c$ and $d$, each of length~$m$, define $\phi(c,d) = b_{1} \cdots b_{m}$, where
\[b_{j}  
%= 1^{c_j+d_j}0
= \begin{cases} 0,& \text{if}\; (c_{j},d_{j}) = (0,0)\\
												10, & \text{if}\; (c_{j},d_{j}) = (1,0) \;\text{or}\; (0,1)\\
												110, & \text{if} \;(c_{j},d_{j}) = (1,1). \end{cases} \]
Note that $\phi$ is a left inverse of $\psi$, i.e.\ $\phi \circ \psi(a) = a$, up to possible concatenation of an extra $0$.											
												
\begin{lemma}
\label{phi}
Let $Y$ be a principal branch of the infinite ternary tree.  Fix a rotor configuration on $Y$ with the root rotor pointing to $o$.  Let $c$ and $d$ be the escape sequences for the configurations on the left and right sub-branches of $Y$, respectively.  Then $\phi(c,d)$ is the escape sequence for the full branch $Y$.
\end{lemma}

\begin{proof}
We claim that each word $b_j$ is the escape sequence for the $j^{th}$ full rotation of the root rotor.   Note that after the root rotor has performed $j-1$ full rotations, each of the sub-branches $L$ and $R$ of $Y$ has seen exactly $j-1$ chips, so the next chip sent to $L$ (resp.\ R) will either return to $r$ or escape to infinity accordingly as $c_j=0$ or $c_j=1$ (resp.\ $d_j=0$ or $d_j=1$).

Consider first the case $(c_j,d_j) = (0,0)$.  After $j-1$ full rotations of the root rotor, the next chip that enters $Y$ will be routed first to $L$, then returned to $r$, sent to $R$, returned to $r$, and finally routed back up to the origin.  The root rotor has now performed a full turn, with corresponding escape sequence $b_j=0$.
If $(c_j,d_j) = (1,0)$, the next chip entering $Y$ will be routed to $L$, where it escapes to infinity. The following chip will be routed to $R$ and then back up to the origin, completing a full rotation of the root rotor.  In this case we have escape sequence $b_j = 10$.  If $(c_j,d_j) = (0,1)$, the next chip entering $Y$ will be routed to $L$, back up to $r$, and then to $R$ where it escapes to infinity. The following chip will be routed directly up to the origin leaving the root rotor pointing up once again.  Again, in this case $b_j = 10$.  Finally, if $(c_j,d_j) = (1,1)$, the next two chips entering $Y$ will escape to infinity, the first through $L$ and the second through $R$. The following chip will be routed directly up to the origin, once again leaving the root rotor pointing up.  In this case we have $b_j = 110$. 
\end{proof}

To adapt Lemma~\ref{phi} to the case when the root rotor is not pointing up, we define \emph{extended escape sequences} $c'$ and $d'$ associated to the two sub-branches.
If the root rotor initially points to $L$, let $c' = 0c$ and $d' = d$.  If the root rotor initially points to $R$, let $c' = 0c$ and $d' = 0d$. Then $a = \phi(c',d')$ is the escape sequence of the full branch $Y$. 
%Note that if the binary words $c$ and $d$ satisfy $(P_{k})$ then so do $c'$ and $d'$.

We now introduce the condition that is central to characterizing which words can be escape sequences:
\renewcommand{\theequation}{$P_k$}
\begin{equation} 
 \text{any subword of length $2^{k}-1$ contains at most $2^{k-1}$ ones}
\end{equation}
We next show that the map $\psi$ preserves this requirement. 
\renewcommand{\theequation}{\arabic{equation}}

\begin{lemma}
\label{psi}
Let $a$ be a binary word satisfying $(P_{k})$ and let $\psi(a) = (c,d)$ as defined in (\ref{psidef}). Then $c$ and $d$ each satisfy $(P_{k-1})$.
\end{lemma}

\begin{proof}
Let $c'$ be a subword of $c$ of length $2^{k-1}-1$ and let $d'$ be the corresponding subword of $d$. Let $a' = \phi(c',d')$, which is a subword of $a0$. The formula for $\phi$ guarantees that $a'$ has one zero for each letter of $c'$, so $a'$ has $2^{k-1}-1$ zeros.  Since the last letter of $a'$ is zero, and $a$ satisfies $(P_k)$, it follows that $a'$ has at most $2^{k-1}$ ones (else after truncating the final zero, the suffix of $a'$ of length $2^k-1$ has at most $2^{k-1}-2$ zeros, hence at least $2^{k-1}+1$ ones).
%Suppose $a'$ has more than $2^{k-1}$ ones. We know the last letter of $a'$ is zero so we remove this zero and look at the subword consisting of the last remaining $2^{k}-1$ letters of $a'$. This subword has at most $2^{k-1}-2$ zeros and therefore more than $2^{k-1}$ ones which contradicts the requirement that $a$ satsify $(P_{k})$. Thus $a'$ has at most $2^{k-1}$ ones.

Let $m$ be the number of ones in $c'$.  Since the instances of $(0,1)$ and $(1,0)$ alternate in the formula for $\psi(a)=(c,d)$, it follows that $d'$ must have at least $m-1$ ones. Since the number of ones in $c'$ and $d'$ combined equals the number of ones in $a'$, we obtain $2m-1 \leq 2^{k-1}$, hence $m \leq 2^{k-2}$.  The same argument with the roles of $c$ and $d$ reversed shows that $d$ has at most $2^{k-2}$ ones.
\end{proof}

\begin{lemma}
\label{branch-escape}
Let $a=a_1\ldots a_n$ be a binary word of length $n$. Then $a$ is an escape sequence for some rotor configuration on the infinite branch $Y$ if and only if $a$ satisfies $(P_{k})$ for all $k$.
\end{lemma}

\begin{proof}
Suppose $a$ is an escape sequence. We prove that $a$ satisfies $(P_{k})$ for each $k$ by induction on $k$. That $a$ satisfies $(P_{1})$ is trivial. Now suppose that every escape sequence satisfies $(P_{k-1})$ and let $c$ and $d$ be the extended escape sequences of the left and right sub-branches respectively.  Then $a = \phi(c,d)$ up to the possible concatenation of an extra zero.  Let $a'$ be a subword of $a$ of length $2^{k}-1$, and let $\psi(a')=(c',d')$.  Then there are words $c''$ and $d''$ each of which is a subword of $c$ or $d$, and which are equal to $c'$ and $d'$, respectively, except possibly in the first letter; moreover the first letters satisfy $c'_1 \leq c''_1$ and $d'_1 \leq d''_1$.

By the formula for $\psi$, the number of ones in $a'$ is the sum of the number of ones in $c'$ and $d'$. If $c'$ has length at most $2^{k-1}-1$, then since $c$ and $d$ satisfy $(P_{k-1})$, each of $c'$ and $d'$ has at most $2^{k-2}$ ones, and therefore $a'$ has at most $2^{k-1}$ ones.  On the other hand, if $c'$ has length at least $2^{k-1}$, then the number of zeros in $a'$ is at least $2^{k-1}-1$. Thus $a'$ has at most $2^{k-1}$ ones, so $a$ satisfies $(P_{k})$.
%to Then $c'$ and $d'$ combined have more than $2^{k-1}$ ones and therefore at least one of these strings has more than $2^{k-2}$ ones. Say it is $c'$ without loss of generality. Since $|c'| - |d'| \le 1$, $c'$ must have at most as many zeros as $d'$ does. Let $z$ be the number of zeros in $a'$, $2*o$ the number of ones in $a'$, and $d$ the difference ($d=2*o-z$). We see from the definition of $\phi$ that $a'$ has at least $2*o-z$ subwords of "`110"' so $c'$ and $d'$ each have at most $z-(2*o-z)=2z-2o$ zeros. Thus $c'$ has at least $o \ge 2^{k-2}$ ones and at most $2z-2o \le 2^{k}-2 - 2^{k-1} = 2^{k-1} - 2$.

The proof of the converse is by induction on $n$. For $n=1$ the statement is trivial. Suppose that every binary word of length $n-1$ satisfying $(P_{k})$ for each $k$ is an escape sequence. Then by Lemma~\ref{psi}, $\psi(a) = (c,d)$ gives a pair of binary words each satisfying $(P_{k})$ for all $k$. If $c$ and $d$ have length $n-1$ or less, then they are escape sequences by induction, hence $a$ is an escape sequence by Lemma~\ref{phi}. If $c$ and $d$ are of length~$n$, then the definition of~$\psi$ implies that $a_j=0$ for all $j<n$, in which case $a$ is an escape sequence by the remark following Proposition~\ref{infinitetreesreturn}.
\end{proof}

We can now establish our main result characterizing all possible escape sequences on the infinite ternary tree.

\begin{theorem}
\label{escape}
Let $a=a_1\ldots a_n$ be a binary word. For $j \in \{1,2,3\}$ write $a^{(j)} = a_{j} a_{j+3} a_{j+6} \ldots$. Then $a$ is an escape sequence for some rotor configuration on the infinite ternary tree $T$ if and only if each $a^{(j)}$ satisfies $(P_{k})$ for all $k$.
\end{theorem}

\begin{proof}
Let $Y^{(1)}$, $Y^{(2)}$, and $Y^{(3)}$ be the three principal branches of $T$ assigned so that the rotor at the origin initially points to $Y^{(3)}$. Then $a$ is the escape sequence for $T$ if and only if $a^{(j)} = a_{j} a_{j+3} a_{j+6} \ldots $ is the escape sequence for $Y^{(j)}$. The result now follows from Lemma~\ref{branch-escape}.
\end{proof}

\section{Concluding Remark}

We conclude with an open question.  While Theorem~\ref{escape} completely characterizes the possible escape sequences for rotor-router walk on the infinite ternary tree, we know nothing about the possible escape sequences for rotor-router walk on another natural class of transient graphs, namely $\Z^d$ for $d\geq 3$.  The open question is this: does there exist a rotor configuration on $\Z^d$ for $d\geq 3$, analogous to the configuration on the tree described in Proposition~\ref{infinitetreesreturn}, so that every chip in an infinite sequence of chips started at the origin eventually returns to the origin?  We remark that Jim Propp has found such a configuration on $\Z^2$.

\section*{Acknowledgments}

The authors thank Yuval Peres, Jim Propp and Parran Vanniasegaram for useful discussions.

\end{document}